\documentclass[12pt, reqno]{amsart}
\usepackage{amsmath, amsthm, amscd, amsfonts, amssymb, graphicx, color, enumerate}
\usepackage[bookmarksnumbered, colorlinks, plainpages]{hyperref}
\usepackage[all,cmtip]{xy}
\usepackage{mathrsfs}
\hypersetup{colorlinks=true,linkcolor=red, anchorcolor=green, citecolor=cyan, urlcolor=red, filecolor=magenta, pdftoolbar=true}

\allowdisplaybreaks[4]

\newtheorem{theorem}{Theorem}[section]

\newtheorem{corollary}[theorem]{Corollary}
\theoremstyle{definition}

\newtheorem{question}[theorem]{Question}

\numberwithin{equation}{section}

\begin{document}
\setcounter{page}{1}

\title[Essential norms of composition operators]{Essential norms of composition operators on Hardy spaces of homogeneous rooted trees}

\author[L. Liu]{Lei Liu$^{*}$}

\address{School of Mathematics and Statistics, Xidian University, Xi'an 710126, P.R. China.}

\email{\textcolor[rgb]{0.00,0.00,0.84}{leiliu@mail.xidian.edu.cn}}

\author[W. Ma]{Wanxin Ma}

\address{School of Mathematics and Statistics, Xidian University, Xi'an 710126, P.R. China.}

\email{\textcolor[rgb]{0.00,0.00,0.84}{mawanxin26@126.com}}

\makeatletter
\@namedef{subjclassname@2020}{\textup{2020} Mathematics Subject Classification}
\makeatother


\date{Received: xxxxxx; Revised: yyyyyy; Accepted: zzzzzz.
\newline \indent $^{*}$ Corresponding author}

\begin{abstract}
We establish characterizations of the essential norms of composition operators on Hardy spaces of homogeneous rooted trees. As direct applications of our main results, we obtain compactness criteria for such composition operators.
\end{abstract} \maketitle

\section{Introduction}

Let $X$ be a Banach space consisting of complex-valued functions on a set $\Omega$. Given a self-map $\varphi:\Omega\to\Omega$, the composition operator $C_{\varphi}:X\to X$ is defined by
\[
C_{\varphi}(f)=f\circ\varphi, \ \ \ \forall f\in X.
\]
Composition operators have been extensively investigated on various Banach spaces of holomorphic functions over the past half century. The main subject in the study of composition operators is to describe operator theoretic properties of $C_{\varphi}$ in terms of function theoretic properties of $\varphi$. For related work, we refer to \cite{chl,cm,fgw,sh} and the references therein.

One central problem in the theory of composition operators is to characterize their essential norms. Recall that the essential norm $\|A\|_e$ of a bounded operator $A:X\to X$ is defined as the distance from $A$ to the space of compact operators on $X$. For example, the essential norms of composition operators between holomorphic function spaces were characterized in \cite{gll,no,sh2,zh}.

In recent years, considerable work has been devoted to the study of composition operators on function spaces over discrete structures such as trees.
Allen, Colonna and Easley \cite{ace} characterized boundedness, compactness, isometries and spectra of composition operators on Lipschitz space of a rooted tree. Muthukumar and Ponnusamy \cite{mp} studied the bounded, compact, invertible, isometric composition operators on Hardy spaces of homogeneous rooted trees. Colonna and  Mart\'{i}nez-Avenda$\tilde{\mathrm{n}}$o investigated boundedness, spectra and hypercyclicity of composition operators on the little Lipschitz space of a rooted tree in the recent paper \cite{cma}.
To the best of our knowledge, the essential norm estimates for composition operators on discrete function spaces have not previously appeared in the literature. Thus we naturally consider the following question.
\begin{question}\label{que}
Is it possible to establish estimates for the essential norms of composition operators acting on function spaces over trees?
\end{question}
In this paper, we answer the above question by deriving explicit formulas for the essential norms of composition operators on Hardy spaces of homogeneous rooted trees.

The paper is organized as follows. We refer to the Section 2 for some preliminaries on trees and the definitions of Hardy spaces of homogeneous rooted trees. In Section 3, we derive formulas for the essential norms of composition operators on Hardy spaces of homogeneous rooted trees. As immediate consequences, we give compactness criteria for composition operators on these Hardy spaces. It is worth noting that Corollary \ref{cor} gives the first necessary and sufficient characterization for the compactness of the composition operators on Hardy spaces with $1\leq p<\infty$.

\section{preliminaries}

A tree is a locally finite connected and simply connected graph $T$, which, as a set, we identify with the collection of its vertices. Two vertices $u$ and $v$ are said to be neighbors if there is an edge connecting them, and we adopt the notation $u \sim v$. The degree of a vertex is the number of its neighbors. A tree is said to be $k$-homogeneous if every vertex has degree exactly $k$. A path between two vertices $u$ and $w$ is a finite sequence of distinct vertices $[u, u_1, u_2, \ldots, u_{k}, w]$ such that $u \sim u_1 \sim u_2 \sim \cdots \sim u_{k} \sim w$. Since $T$ is a tree, the path between two vertices is unique.
	
Given a tree $T$ rooted at $o$ and a vertex $u\in T$, a vertex $v \neq u$ is called a descendant of $u$ if $u$ belongs to the unique path from $o$ to $v$. In this case, $u$ is called an ancestor of $v$. We write $v^-$ for the parent of $v$, namely the unique neighbor of $v$ which is an ancestor of $v$. The vertex $v$ is called a child of $v^-$. The length of a path between two vertices $u$ and $w$ is the number of edges connecting them. The distance $d(u, w)$ between vertices $u$ and $w$ is the length of the path connecting $u$ to $w$. In particular, the length of a path joining the root $o$ to a vertex $w$ is called the length of $w$ and is denoted by $|w|$.	
A function on a tree is a complex-valued function on the set of its vertices.

Throughout this paper, $T$ always stands for a $(q+1)$-homogeneous rooted tree with root vertex $o$, and we adopt the convention $\mathbb{N}_{0}=\mathbb{N} \cup \{0\}$. For each $n\in\mathbb{N}$, we define
\[
M_{p}(n, f):= \begin{cases}\left(\frac{1}{(q+1) q^{n-1}} \sum_{|v|=n}|f(v)|^{p}\right)^{\frac{1}{p}} & \text{if } p \in(0, \infty), \\ \max _{|v|=n}|f(v)| & \text { if } p=\infty \end{cases}
\]
and
\[
\| f\| _{p}:=\sup _{n \in \mathbb{N}_{0}} M_{p}(n, f),
\]
with the convention $M_{p}(0, f):=|f(o)|$.
The discrete Hardy space associated with the tree $T$, denoted by $\mathbb{T}_{p}$, is defined as
\[
\mathbb{T}_{p}:=\left\{f: T \to \mathbb{C} \mid\| f\| _{p}<\infty\right\}
\]
for all $p \in(0, \infty]$. Similarly, the discrete little Hardy space, denoted by $\mathbb{T}_{p,0}$, is defined by
\[
\mathbb{T}_{p,0}:=\left\{f\in\mathbb{T}_{p} \mid \lim_{n\to\infty}M_{p}(n, f)=0\right\}
\]
By \cite[Theorem A]{mp}, both $\mathbb{T}_{p}$ and $\mathbb{T}_{p,0}$ form Banach spaces for $1\leq p\leq\infty$. For the sake of simplicity, we omit the parameter $q$ from the notation for these spaces.
	
We fix the following notation for the remainder of the paper. $\|\cdot\|$ denotes $\|\cdot\|_p$ on $\mathbb{T}_{p}$ spaces. For each $n \in \mathbb{N}_{0}$, write $D_{n}$ for the collection of vertices $v \in T$ satisfying $|v|=n$, and let $c_{n}$ denote the cardinality of $D_{n}$. Accordingly,
\[
c_{n}= \begin{cases}(q+1) q^{n-1} & \text { if } n \in \mathbb{N}, \\ 1 & \text { if } n=0. \end{cases}
\]
For $n\in\mathbb{N}_{0}$ and $w\in T$, let $N_{\varphi}(n,w)$ stand for the number of preimages of $w$ under $\varphi$ lying in the set $D_{n} = \{v\in T \mid |v|=n\}$. Equivalently,
\[
N_{\varphi}(n,w) = \#\big(\varphi^{-1}(w)\cap D_{n}\big).
\]
Furthermore, for all $m,n\in\mathbb{N}_{0}$, define $N_{m,n}$ as the maximum of $N_{\varphi}(n,w)$ taken over all vertices $w$ with $|w|=m$. It is straightforward to verify that
\[
\sum_{m=0}^{\infty} N_{m,n} \leq c_{n}
\]
holds for every nonnegative integer $n$.

\section{Results and proofs}

In this section, we discuss essential norms of composition operator $C_{\varphi}$ on $\mathbb{T}_{p}$ spaces. We state the first main result as follows.

\begin{theorem}\label{thm1}
Let $T$ be a homogeneous rooted tree, and let $1\leq p<\infty$. Suppose that $\varphi$ is a self-map of $T$ such that the composition operator $C_{\varphi}: \mathbb T_p\to\mathbb T_p$ is bounded. Then the essential norm of $C_{\varphi}$ is
\[
\|C_{\varphi}\|_e=\lim_{M\to\infty}\{\sup_{n\in\mathbb N_0}\frac1{c_n}\sum_{m=M}^{\infty}N_{m,n}c_m\}^{1/p}.
\]
\end{theorem}
	
\begin{proof}
The upper estimate.
For any $M\in\mathbb N$, define $P_M: \mathbb{T}_p\to\mathbb{T}_p$ by
		\[
		(P_Mf)(v)=
		\begin{cases}
			f(v), & |v|<M,\\
			0, & |v|\ge M.
		\end{cases}
		\]
This yields that $P_M$ is a finite-rank operator. Put $K_M=C_{\varphi}P_M$. Then $K_M$ is compact, and hence
		\begin{equation}
		\|C_{\varphi}\|_e
		\le \|C_{\varphi}-K_M\|
		=\|C_{\varphi}(I-P_M)\|.
		\end{equation}
For any $f\in \mathbb{T}_p$ and $n\in\mathbb N_0$, we have
		\begin{align*}
			M_p(n,C_{\varphi}(I-P_M)f)^p
			&=
			\frac1{c_n}\sum_{|v|=n}
			|(I-P_M)f(\varphi(v))|^p \\
			&=
			\frac1{c_n}\sum_{m=M}^{\infty}\sum_{|w|=m}
			N_{\varphi}(n,w)|f(w)|^p \\
			&\le
			\frac1{c_n}\sum_{m=M}^{\infty}N_{m,n}
			\sum_{|w|=m}|f(w)|^p \\
			&\le
			\frac1{c_n}\sum_{m=M}^{\infty}N_{m,n}c_m\,\|f\|_{p}^p .
		\end{align*}
It follows that
		\begin{equation}
		\|C_{\varphi}(I-P_M)\|^p
		\le
		\sup_{n\in\mathbb N_0}
		\frac1{c_n}\sum_{m=M}^{\infty}N_{m,n}c_m .
		\end{equation}
Combining (3.1) and (3.2), we obtain that
		\[
		\|C_{\varphi}\|_e
		\le
		\lim_{M\to\infty}\{\sup_{n\in\mathbb N_0}\frac1{c_n}\sum_{m=M}^{\infty}N_{m,n}c_m\}^{1/p}.
		\]
		
The lower estimate. For any $M\in\mathbb N$, we may choose $n_M\in\mathbb N_0$ such that
		\begin{equation}
		\frac1{c_{n_M}}\sum_{m=M}^{\infty}N_{m,n_M}c_m>\sup_{n\in\mathbb N_0}
		\frac1{c_n}\sum_{m=M}^{\infty}N_{m,n}c_m-\frac1M .
		\end{equation}
For each $m\ge M$ with $N_{m,n_M}>0$, pick $w_{M,m}\in D_m$ satisfying
		\[
		N_{\varphi}(n_M,w_{M,m})=N_{m,n_M}.
		\]
Define
		\[
		g_M :=
		\sum_{\substack{m\ge M\\N_{m,n_M}>0}}
		c_m^{1/p}\chi_{\{w_{M,m}\}} .
		\]
For any $M\in\mathbb N$, since
		\[
		\sum_{m=0}^{\infty}N_{m,n_M}\le c_{n_M},
		\]
there are only finitely many $m$ for which $N_{m,n_M}>0$. Hence $g_M\in\mathbb{T}_{p,0}\subset\mathbb{T}_p$.
Moreover, for any $m\in\mathbb N_0$, we have
		\[
		M_p(m,g_M)=
		\begin{cases}
			1, & m\ge M\text{ and }N_{m,n_M}>0,\\
			0, & \text{otherwise},
		\end{cases}
		\]
witch implies $\|g_M\|_{p}\le1$.
		
Now we show that $g_M\overset{w}{\rightarrow}0$ as $M\to\infty$. Indeed, for each $m\in\mathbb N_0$, let
		\[
		E_m=\{h:D_m\to\mathbb C\}
		\]
with
		\[
		\|h\|_{E_m}=
		\left(
		\frac1{c_m}\sum_{|v|=m}|f(v)|^p
		\right)^{1/p}.
		\]
Then, for any $f\in\mathbb{T}_{p,0}$, we have
		\[
		\|f|_{D_m}\|_{E_m}=M_p(m,f)\to0
		\]
as $M\to\infty$	and
		\[
		\|f\|_{p}=\sup_{m\in\mathbb N_0}\|f|_{D_m}\|_{E_m}.
		\]
Therefore,
		\[
		\mathbb{T}_{p,0}\cong
\bigoplus\nolimits_{0}E_m:=\left\{x\in\prod_{m=0}^{\infty}E_m:\|x(m)\|_{E_m}\to0\right\}.
		\]
By the standard duality for direct sums, we have
		\begin{equation}
		\left(
		\bigoplus\nolimits_{0}E_m
		\right)^{*}
		\cong
		\bigoplus\nolimits_{1}E_m^{*}.
		\end{equation}
For any $F\in(\mathbb{T}_p)^*$ and $f\in\mathbb{T}_{p,0}\subset\mathbb{T}_{p}$, it follows from (3.4) that there exist $F_m\in E_m^*$, $m\in\mathbb N_0$ such that		
		\[
		F(f)
		=
		\sum_{m=0}^{\infty}F_m(f|_{D_m}),
		\qquad
		\sum_{m=0}^{\infty}\|F_m\|<\infty.
		\]
		Since $g_M|_{D_m}=0$ for $m<M$, we obtain that
		\begin{align*}
			|F(g_M)|
			&\le
			\sum_{m=M}^{\infty}
			|F_m(g_M|_{D_m})| \\
			&\le
			\sum_{m=M}^{\infty}
			\|F_m\|\,
			\|g_M|_{D_m}\|_{E_m} \\
			&\le
			\sum_{m=M}^{\infty}\|F_m\|.
		\end{align*}
Since $\sum_{m=0}^{\infty}\|F_m\|<\infty$, it follows that $F(g_M)\to0$ as $M\to\infty$. This yields
		\[
		g_M\overset{w}{\rightarrow}0
		\]
as $M\to\infty$ in $\mathbb{T}_p$.
		
Furthermore, by the definition of $g_M$ and the choice of $w_{M,m}$, we have
		\begin{align*}
			M_p^p(n_M,C_{\varphi}g_M)
			&=
			\frac1{c_{n_M}}
			\sum_{|v|=n_M}|g_M(\varphi(v))|^p \\
			&=
			\frac1{c_{n_M}}
			\sum_{\substack{m\ge M\\N_{m,n_M}>0}}
			c_mN_{\varphi}(n_M,w_{M,m}) \\
			&=
			\frac1{c_{n_M}}
			\sum_{\substack{m\ge M\\N_{m,n_M}>0}}
			N_{m,n_M}c_m \\
			&=
			\frac1{c_{n_M}}
			\sum_{m=M}^{\infty}N_{m,n_M}c_m.
		\end{align*}
Therefore, by (3.3), we arrive at
		\[
		\|C_{\varphi}g_M\|_{p}^p
		\ge
		\frac1{c_{n_M}}
		\sum_{m=M}^{\infty}N_{m,n_M}c_m
		>
		\sup_{n\in\mathbb N_0}
		\frac1{c_n}\sum_{m=M}^{\infty}N_{m,n}c_m
		-\frac1M .
		\]
For any compact operator $K$, since $g_M\overset{w}{\rightarrow}0$ in $\mathbb{T}_p$, it follows that
		\[
		\|Kg_M\|_{p}\to0
		\]
as $M\to\infty$. Thus,
		\begin{align*}
			\|C_{\varphi}-K\|
			&\ge
			\limsup_{M\to\infty}
			\|(C_{\varphi}-K)g_M\|_{p} \\
			&\ge
			\limsup_{M\to\infty}
			\bigl(\|C_{\varphi}g_M\|_{p}-\|Kg_M\|_{p}\bigr) \\
			&=
			\limsup_{M\to\infty}\|C_{\varphi}g_M\|_{p}.
		\end{align*}
Taking $p$-th powers, we obtain that
		\begin{align*}
			\|C_{\varphi}-K\|^p
			&\ge
			\limsup_{M\to\infty}\|C_{\varphi}g_M\|_{p}^p \\
			&\ge
			\limsup_{M\to\infty}
			\{
			\sup_{n\in\mathbb N_0}
			\frac1{c_n}\sum_{m=M}^{\infty}N_{m,n}c_m
			-\frac1M
			\} \\
			&=
			\lim_{M\to\infty}
			\{\sup_{n\in\mathbb N_0}
			\frac1{c_n}\sum_{m=M}^{\infty}N_{m,n}c_m\}.
		\end{align*}
Since $K$ is an arbitrary compact operator, we have
		\[
		\|C_{\varphi}\|_e
		\ge
		\lim_{M\to\infty}
		\{\sup_{n\in\mathbb N_0}
		\frac1{c_n}\sum_{m=M}^{\infty}N_{m,n}c_m\}^{1/p}.
		\]
Combining the lower and upper estimate, we arrive at
		\[
		\|C_{\varphi}\|_e
		=
		\lim_{M\to\infty}
		\{\sup_{n\in\mathbb N_0}
		\frac1{c_n}\sum_{m=M}^{\infty}N_{m,n}c_m\}^{1/p}.
		\]
The proof is complete.
	\end{proof}

As an immediate consequence of Theorem \ref{thm1}, we obtain the following corollary, which characterizes the compactness of composition operators on $\mathbb{T}_p$.

\begin{corollary}\label{cor}
Let $T$ be a homogeneous rooted tree, and let $1\leq p<\infty$. Suppose that $\varphi$ is a self-map of $T$ such that $C_{\varphi}$ is bounded on $\mathbb{T}_p$. Then the following conditions are equivalent:
\begin{enumerate}[(i)]
\item $C_{\varphi}$ is compact on $\mathbb{T}_p$;
\item $\lim_{M\to\infty}\{\sup_{n\in\mathbb{N}_0}
\frac{1}{c_n}\sum_{m=M}^{\infty}N_{m,n}c_m\}=0$.
\end{enumerate}
\end{corollary}

Now we consider the essential norms of composition operators on $\mathbb T_p$ for $1\leq p<\infty$ over $2$-homogeneous trees. Note that every self-map $\varphi$ of $2$-homogeneous trees induces a bounded $C_{\varphi}$ on $\mathbb T_p$ (see \cite[Theorem 2]{mp2}).

\begin{theorem}\label{thm2}
Let $T$ be a $2$-homogeneous tree, and let $D_n=\{a_n,b_n\}$ for each $n\in\mathbb N$. Suppose that $\varphi$ is a self-map of $T$ which induces the composition operator $C_{\varphi}: \mathbb T_p\to\mathbb T_p$ for $1\le p<\infty$. Then we have the following:
\begin{enumerate}[(i)]
\item If $\varphi$ is bounded, then $\|C_{\varphi}\|_e=0$.
\item If $\varphi$ is unbounded, then one of the following case occurs:
\begin{enumerate}[(a)]
\item If there exists $N\in\mathbb N$ such that for every $n\in\mathbb N$, whenever $|\varphi(a_n)|\ge N$ and $|\varphi(b_n)|\ge N$, $\varphi$ maps $D_n$ onto $D_k$ for some $k\ge N$, then $\|C_{\varphi}\|_e=1$.
\item If for every $M\in\mathbb N$, there exists $n\in\mathbb N$ such that $|\varphi(a_n)|\ge M$ and $|\varphi(b_n)|\ge M$, and either $\varphi(a_n)=\varphi(b_n)$ or $|\varphi(a_n)|\ne|\varphi(b_n)|$, then
$\|C_{\varphi}\|_e=2^{1/p}$.
\end{enumerate}
\end{enumerate}
\end{theorem}

\begin{proof}
For a $2$-homogeneous tree, since $c_n=2$ for all $n\ge1$, Theorem \ref{thm1} implies that
\begin{equation}
\|C_{\varphi}\|_e=\lim_{M\to\infty}\{\sup_{n\in\mathbb N_0}
\sum_{m=M}^{\infty}N_{m,n}\}^{1/p}.
\end{equation}

Suppose that $\varphi$ is bounded. Then there exists $N\in\mathbb N$ such that $|\varphi(v)|<N$ for all $v\in T$. Hence, for any $M\ge N$, we obtain $N_{m,n}=0$ whenever $m\ge M$ and $n\in\mathbb N$. It follows that $\sup_{n\in\mathbb N}\sum_{m=M}^{\infty}N_{m,n}=0$ for all $M\ge N$, which implies $\|C_{\varphi}\|_e=0$ by (3.5).

Next, suppose that $\varphi$ is unbounded and condition $(a)$ holds. We first fix $n\in\mathbb N$ and $M\ge\max\{N,|\varphi(o)|+1\}$. If $|\varphi(a_n)|<M$ and $|\varphi(b_n)|<M$, then $\sum_{m=M}^{\infty}N_{m,n}=0$. Moreover, if $|\varphi(a_n)|<M$ and $|\varphi(b_n)|\geq M$, or $|\varphi(a_n)|\geq M$ and $|\varphi(b_n)|<M$, then $\sum_{m=M}^{\infty}N_{m,n}=1$. Lastly, if $|\varphi(a_n)|\ge M$ and $|\varphi(b_n)|\ge M$, then the assumption implies that $\varphi$ maps $D_n$ onto $D_k$ for some $k\ge M$. Hence $N_{k,n}=1$ and $N_{m,n}=0$ for all $m\ge M$, $m\ne k$, so $\sum_{m=M}^{\infty}N_{m,n}=1$. Therefore
\[
\sup_{n\in\mathbb N_0}\sum_{m=M}^{\infty}N_{m,n}\le1.
\]
On the other hand, note that $\varphi$ is unbounded. For any $M\ge \max\{N,|\varphi(o)|+1\}$, there exists $v\in T$ such that $|\varphi(v)|\ge M$, which implies
\[
\sup_{n\in\mathbb N_0}\sum_{m=M}^{\infty}N_{m,n}\geq1.
\]
It follows that
\[
\sup_{n\in\mathbb N_0}
\sum_{m=M}^{\infty}N_{m,n}=1
\]
for all $M\ge \max\{N,|\varphi(o)|+1\}$. Hence, by equation (3.5), we have $\|C_\varphi\|_e=1$.

Finally, suppose that $\varphi$ is unbounded and condition $(b)$ holds. Fix $M\in\mathbb N$ and choose $n\in\mathbb N$ as in the assumption. If $\varphi(a_n)=\varphi(b_n)=w$, then $|w|\ge M$ and $N_{|w|,n}=2$, which yields $\sum_{m=M}^{\infty}N_{m,n}=2$. If $|\varphi(a_n)|\ne|\varphi(b_n)|$, write $|\varphi(a_n)|=m_1$ and $|\varphi(b_n)|=m_2$, where $m_1,m_2\ge M$ and $m_1\ne m_2$. Then $N_{m_1,n}=N_{m_2,n}=1$, which implies $\sum_{m=M}^{\infty}N_{m,n}=2$. Therefore
\[
\sup_{n\in\mathbb N_0}\sum_{m=M}^{\infty}N_{m,n}\ge2.
\]
On the other hand, $\sum_{m=0}^{\infty}N_{m,n}\le c_n=2$ for all $n\in\mathbb N$. It follows that
\[
\sup_{n\in\mathbb N_0}\sum_{m=M}^{\infty}N_{m,n}=2
\]
for all $M\in\mathbb N$. Using equation (3.5) again, we see that $\|C_{\varphi}\|_e=2^{1/p}$.
This completes the proof.
\end{proof}

As an application of Theorem \ref{thm2}, we immediately obtain the following equivalent conditions for the composition operator $C_{\varphi}$ to be compact on the $\mathbb T_p$ spaces over $2$-homogeneous trees.
It is exactly the result in \cite[Theorem F]{mp}.

\begin{corollary}
Let $T$ be a $2$-homogeneous tree, and let $\varphi$ be a self-map of $T$. Then composition operator $C_\varphi$ is compact on $\mathbb T_p$ if and only if $\varphi$ is bounded.
\end{corollary}

Finally, we investigate essential norms of composition operators on the $\mathbb T_\infty$ spaces.

\begin{theorem}\label{thm3}
Let $T$ be a homogeneous rooted tree. Suppose that $\varphi$ is a self-map of $T$ which induces the composition operator $C_{\varphi}: \mathbb T_\infty\to\mathbb T_\infty$. Then one of the following holds:
\begin{enumerate}[(i)]
\item $\|C_{\varphi}\|_e=0$ whenever $\varphi$ is bounded;
\item $\|C_{\varphi}\|_e=1$ whenever $\varphi$ is unbounded.
\end{enumerate}
\end{theorem}

\begin{proof}
By \cite[Theorem D]{mp}, every self-map $\varphi$ of $T$ induces a bounded composition operator on $\mathbb T_\infty$ with $\|C_\varphi\|=1$. Suppose that $\varphi$ is bounded. Let $\varphi(T)=\{w_1,\ldots,w_r\}$.
It follows that
\[
C_{\varphi}f
=
\sum_{j=1}^{r}f(w_j)\chi_{\varphi^{-1}(w_j)}
\]
for all $f\in\mathbb{T}_{\infty}$. This implies that $C_{\varphi}$ is finite-rank, and hence $\|C_{\varphi}\|_e=0$.

Now suppose that $\varphi$ is unbounded. Since $\|C_{\varphi}\|=1$, we have $\|C_{\varphi}\|_e\le1$. Thus it suffices to show that $\|C_{\varphi}\|_e\ge1$. Indeed, choose distinct vertices $w_j\in\varphi(T)$ and put $f_j=\chi_{\{w_j\}}$. Then $\|f_j\|_{\infty}=1$ for all $j$. Let $K$ be an arbitrary compact operator on $\mathbb{T}_{\infty}$. There exists a subsequence $\{f_{j_k}\}_{k=1}^{\infty}$ such that $\{Kf_{j_k}\}_{k=1}^{\infty}$ converges in norm.
For any $k\in\mathbb{N}$, set $g_k=f_{j_{2k}}-f_{j_{2k+1}}$. Then $\|g_k\|_{\infty}=1$ and
\[
\|Kg_k\|_{\infty}
=
\|Kf_{j_{2k}}-Kf_{j_{2k+1}}\|_{\infty}
\to 0
\]
as $k\to\infty$.

Moreover, for each $k\in\mathbb{N}$, choose $v_k\in T$ such that $\varphi(v_k)=w_{j_{2k}}$. Since $w_{j_{2k}}\ne w_{j_{2k+1}}$, we have $(C_{\varphi}g_k)(v_k)=1$, and hence $\|C_{\varphi}g_k\|_{\infty}=1$. Therefore
\begin{align*}
\|C_{\varphi}-K\|
&\ge \|(C_{\varphi}-K)g_k\|_{\infty}\\
&\ge 1-\|Kg_k\|_{\infty}.
\end{align*}
Letting $k\to\infty$, we obtain $\|C_{\varphi}-K\|\ge1$. Since $K$ is an arbitrary compact operator, it follows that $\|C_{\varphi}\|_e\ge1$. The proof is complete.
\end{proof}

Using Theorem \ref{thm2}, we immediately obtain the following corollary, which restates \cite[Theorem E]{mp}.

\begin{corollary}
Let $T$ be a homogeneous rooted tree. Let $\varphi$ be a self-map of $T$. Then composition operator $C_\varphi$ is compact on $\mathbb T_\infty$ if and only if $\varphi$ is bounded.
\end{corollary}

{\bf Acknowledgement} This work is supported by the Natural Science Basic Research Program of Shaanxi Province (Program No. 2025JC-YBMS-024).

{\bf Author Contributions} All authors share credit for this work equally.

{\bf Data Availability} No datasets were generated or analysed during the current study.

{\bf Declarations}

{\bf Conflict of interest} The authors declare no competing interests.

\bibliographystyle{amsplain}

\end{document}